\documentclass[10pt]{article}

\usepackage{amsthm}
\usepackage{amssymb}
\usepackage{amsmath}
\usepackage{hhline}
\usepackage[mathscr]{eucal}

\ifx\pdfoutput\undefined
   \usepackage[dvips]{graphicx}
   \else
   \usepackage[pdftex]{graphicx}
   \fi
   \usepackage{epstopdf}

\theoremstyle{plain}
\newtheorem{tetel}{Theorem}[section]
\newtheorem{all}[tetel]{Proposition}
\newtheorem{lemma}[tetel]{Lemma}

\newtheorem{sejt}[tetel]{Conjecture}
\theoremstyle{definition}\newtheorem{Def}[tetel]{Definition}
\theoremstyle{remark}\newtheorem{megj}[tetel]{Remark}
\newtheorem{pelda}[tetel]{Example}

\newcommand*{\Z}{\ensuremath{\mathbf Z}}
\newcommand*{\R}{\ensuremath{\mathbf R}}
\newcommand*{\di}{\mathrm d}

\DeclareMathOperator{\id}{id}

\begin{document}

\DeclareGraphicsExtensions{.eps,.pdf}

\title{Braid-positive Legendrian links}
\author{Tam\'as K\'alm\'an}
\maketitle

%
%

\begin{abstract}
Any link that is the closure of a positive braid has a
natural Legendrian representative. These were introduced in \cite{en},
where their Chekanov--Eliashberg contact homology was also evaluated. In
this paper we re-phrase and improve that computation using a matrix
representation. In particular, we present a way
of finding all augmentations of such Legendrians, construct an augmentation which is also
a ruling, and find surprising links to $LU$--decompositions and Gr\"obner bases.
\end{abstract}

\section{Introduction}

I came across a certain set of Legendrian links while searching for
examples to illustrate the main theorem of my thesis \cite{en}, and they
served that purpose very well. Since then I kept returning to them because
I could always discover something pretty. This paper is a collection of
those findings.

The links in question (see Figure \ref{fig:lagpic}), that I call
\emph{Legendrian closures of positive braids}, denote by $L_\beta$, and represent by
front diagrams $f_\beta$, are
Legendrian representatives of braid-positive links, i.e.\ link types that
can be obtained as the closure of a positive braid $\beta$. (These are not
to be confused with the more general notion of positive link, i.e.\ link
types that can be represented with diagrams whose geometric and algebraic
crossing numbers agree.) In fact I conjecture that $L_\beta$ is essentially
the only Legendrian representative of such a link type, in the following
sense.

\begin{sejt} Any braid-positive Legendrian link is a stabilization of the
corresponding Legendrian closure shown in Figure \ref{fig:lagpic}. In
particular, braid-positive links are Legendrian simple. \end{sejt}

This paper however is not about compiling evidence for this conjecture. Let
us only mention that Etnyre and Honda \cite{EH1} proved it for positive
torus knots, that the set of links treated by Ding and Geiges \cite{geig}
includes many two-component braid-positive links (for example, positive $(2k,2)$ torus links), and that Chekanov's example
\cite{chek} of a non-Legendrian simple knot type is $5_2$, which is the
smallest positive, but not braid-positive knot. Also, by Rutherford's work
\cite{rulpoly}, the Thurston--Bennequin number of $L_\beta$ is maximal in its
smooth isotopy class because the front diagram $f_\beta$ is easily seen to admit
rulings (see section \ref{sec:rul}). Because some (actually, all) of those rulings
are $2$--graded, the maximal Thurston--Bennequin number is only attained along with
rotation number $0$.

Instead, we will concentrate on Legendrian isotopy invariants of
$L_\beta$.  Some of these have been evaluated in \cite{en}, of which the present paper
is a continuation. It is thus assumed that the reader is familiar with sections 2 (basic
notions) and 6 (Legendrian closures of positive braids and their relative contact homology)
of \cite{en}. In this paper, we will re-formulate some of those computations, and get new
results also, by using what we call the path matrix of a positive braid. This construction
is very similar to that of Lindstr\"om \cite{lind}
and Gessel--Viennot \cite{gv}, which is also included in the volume
\cite{book}.

The paper is organized as follows. We review some results of \cite{en} in section
\ref{sec:prelim}, and discuss elementary properties of the path matrix in section
\ref{sec:matrix}.
Then, the main results are Theorem \ref{thm:ideal}, where we compute a new generating
set for the (abelianized) image $I$ of the contact homology differential and its
consequence, Theorem \ref{thm:gauss}, which gives a quick test to decide whether a given
set of crossings is an augmentation. The latter is in terms of an $LU$--decomposition,
i.e.\ Gaussian elimination of the path matrix. In Theorem \ref{thm:grob}, we point out
that the old generators, i.e. the ones read off of the knot diagram, automatically form
a Gr\"obner basis for $I$. In Theorem \ref{thm:szimultan}, we construct a subset of the
crossings of $\beta$ which is simultaneously an augmentation and a ruling of $L_\beta$. This
strengthens the well known relationship between augmentations and rulings. We close the
paper with a few examples.



Acknowledgements: Part of the research was carried out in the summer of 2005 when I visited
the Alfr\'ed
R\'enyi Mathematical Institute in Budapest. It is a great pleasure to thank the Institute
for their hospitality and Andr\'as N\'emethi, Endre Szab\'o, and Andr\'as Sz\H{u}cs for
stimulating discussions. I am grateful to Alexander Stoimenow for providing me with a list
of braid-positive knots up to $16$ crossings. My conversations with Mikhail Kotchetov were
invaluable to the discovery of Theorem \ref{thm:grob}. Last but not least, many thanks to
Supap Kirtsaeng for writing the computer implementation of Theorem \ref{thm:gauss}.

\section{Preliminaries}\label{sec:prelim}

The goal of this section is to recall some results from \cite{en} relevant
to this paper. We will work in the standard contact $3$--space
$\R^3_{xyz}$ with the kernel field of the $1$--form $\di z - y \di x$. We
will use the basic notions of Legendrian knot, Legendrian isotopy, front
($xz$) diagram, Maslov potential, Lagrangian ($xy$) diagram, resolution
\cite{computable}, Thurston--Bennequin ($tb$) and rotation ($r$) numbers,
admissible disc and contact homology\footnote{Because absolute contact
homology doesn't appear in the paper, we'll use this shorter term for what
may be better known as relative, Legendrian, or Chekanov--Eliashberg
contact homology.} etc.\ without reviewing their definitions. We will also
assume that the reader is familiar with section 6 of \cite{en}, of which
this paper is in a sense an extension. For a complete introduction to
Legendrian knots and their contact homology, see \cite{etn}.

We would like to stress a few points only whose treatment may be somewhat
non-standard. Crossings $a$ of both front and Lagrangian diagrams are assigned
an \emph{index}, denoted by $|a|$, which is an element of $\Z_{2r}$, with the entire
assignment known as a \emph{grading}. This is easiest to define for fronts
of single-component knots as the difference of the Maslov potentials
($\text{upper}-\text{lower}$) of the two intersecting strands. If a
Lagrangian diagram is the result of resolution, the old crossings keep their indices and
the crossings replacing the right cusps are assigned
the index $1$. In the multi-component case, the Maslov potential difference
becomes ambiguous for crossings between different components. This gives
rise to an infinite set of so-called admissible gradings. We consider
these as introduced in \cite[section 2.5]{computable} and not the larger class of
gradings described in \cite[section 9.1]{chek}.

\begin{figure}
\centering
\includegraphics[width=\linewidth]{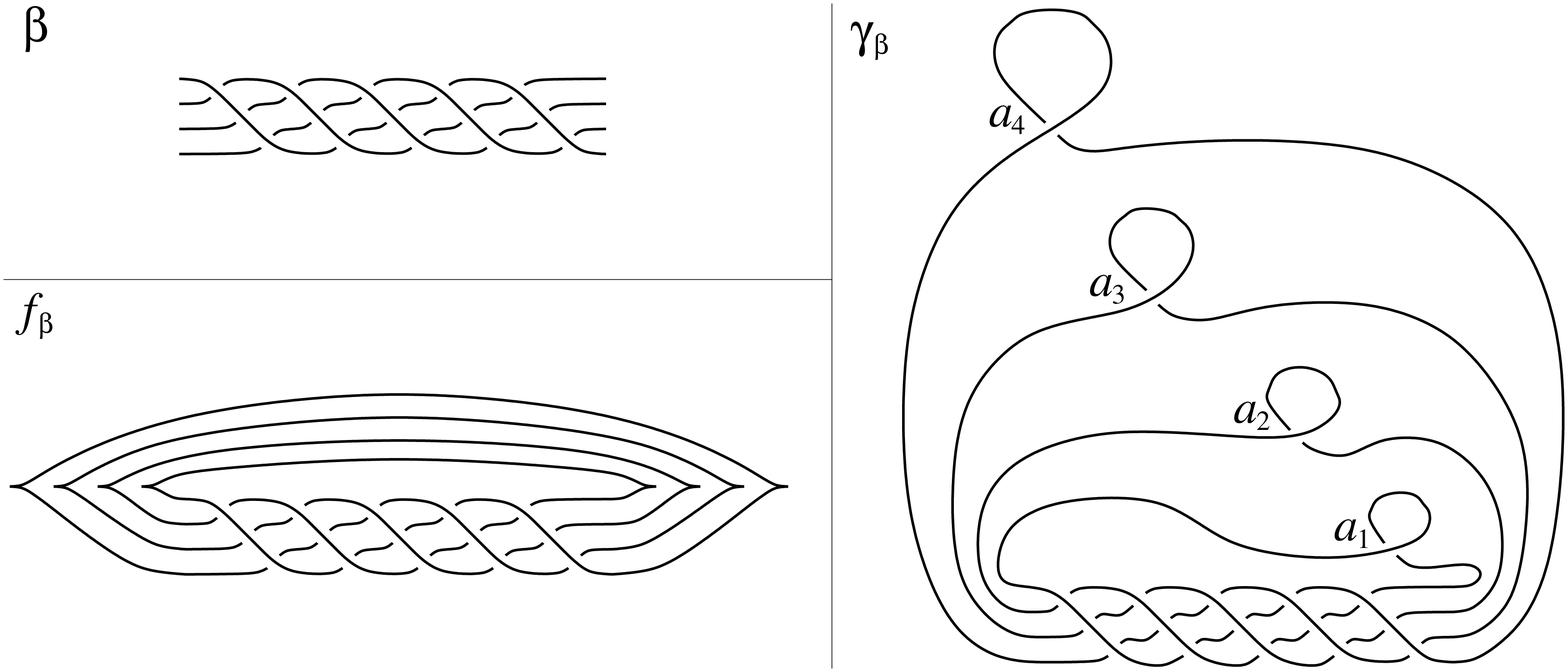}
\caption{Front ($f_\beta$) and Lagrangian ($\gamma_\beta$) diagrams of the
closure ($L_\beta$) of the positive braid $\beta$}
\label{fig:lagpic}
\end{figure}

Let $\beta$ denote an arbitrary positive braid word. The Legendrian
isotopy class $L_\beta$ is a natural Legendrian representative of the link
which is the closure of $\beta$. (All braids and braid words in this paper
are positive. The same symbol $\beta$ may sometimes refer to the braid
represented by the braid word $\beta$.) $L_\beta$, in turn, is represented by the
front diagram
$f_\beta$ and its resolution, the Lagrangian diagram $\gamma_\beta$ (see
Figure \ref{fig:lagpic}). Considering $\beta$ drawn horizontally, label
the left and right endpoints of the strands from top to bottom with the
first $q$ whole numbers ($q$ is the number of strands in $\beta$). The
crossings of $\beta$, labeled from left to right by the symbols
$b_1,\ldots,b_w$, are the only crossings of $f_\beta$. Due to resolution,
$\gamma_\beta$ also has the crossings $a_1,\ldots,a_q$.

The differential graded algebra (DGA) $\mathscr A$ which is the chain
complex for the contact homology of $L_\beta$ is generated (as a
non-commutative algebra with unit) freely over $\Z_2$ by these $q+w$
symbols ($w$ is the word length or exponent sum of $\beta$). It's assigned
a $\Z$--grading\footnote{All components of $L_\beta$ have $r=0$. If there
are multiple components, what we describe here is only one of the
admissible gradings.} which takes the value $0$ on the $b_k$ and the value
$1$ on the $a_n$ (extended by the rule $|uv|=|u|+|v|$). By Theorem
6.7 of \cite{en}, the differential $\partial$ is given on the generators by the
formulas \begin{equation}\label{eq:perem}
\partial(b_k)=0\quad\text{and}\quad\partial(a_n)=1+C_{n,n}. \end{equation}
(It is extended to $\mathscr A$ by linearity and the Leibniz rule.) Here, for any $n$,
\begin{equation}\label{eq:Cii} C_{n,n}
=\sum_{\{\,i_1,\ldots,i_c\,\}\in D_n}
B_{n,i_1}B_{i_1,i_2}B_{i_2,i_3}\ldots B_{i_{c-1},i_c}B_{i_c,n},
\end{equation}
where two more terms require explanation.

\begin{Def}\label{def:sorozat}
A finite sequence of positive integers is called \emph{admissible}
if for all $s\ge 1$, between any two appearances of $s$ in the
sequence there is a number greater than $s$ which appears between
them. For $n\ge 1$, we denote by $D_n$ the set of all admissible
sequences that are composed of the numbers $1,2,\ldots,n-1$.
\end{Def}

Note that non-empty admissible sequences have a unique highest
element.

\begin{Def}\label{def:B} Let $1\le i,j\le q$. The element $B_{i,j}$ of the
DGA of $\gamma_\beta$ is the sum of the following products. For each path
composed of parts of the strands of the braid (word) $\beta$ that connects
the left endpoint labeled $i$ to the right endpoint labeled $j$ so that it
only turns around quadrants facing up, take the product of the labels of
the crossings from left to right that it turns at. (We will refer to the
paths contributing to $B_{i,j}$ as \emph{paths in the braid}.) \end{Def}

We will also use the following notation: for any $i<j$, let
\begin{equation}\label{eq:Cij}
C_{i,j}=\sum_{\{\,i,i_1,\ldots,i_c\,\}\in D_j} B_{i,i_1}B_{i_1,i_2}B_{i_2,i_3}\ldots
B_{i_{c-1},i_c}B_{i_c,j}. \end{equation}

The expressions $B_{i,j}$ and $C_{i,j}$ are elements of the DGA
$\mathscr A$. Even though $\mathscr A$ is non-commutative, we will
refer to them, as well as to similar expressions and even matrices
with such entries, as polynomials.

\section{The path matrix}\label{sec:matrix}

%
%

The polynomials $B_{i,j}$ are naturally arranged in a $q\times q$
matrix $B_\beta$ (with entries in $\mathscr A$), which we will call
the \emph{path matrix of $\beta$}.

If we substitute $0$ for each crossing label of $\beta$, then $B_\beta$
reduces to the matrix of the underlying permutation $\pi$ of $\beta$:
\[B_\beta(0,0,\ldots,0)=\left[\delta_{\pi(i),j}\right]=:P_\pi,\] where
$\delta$ is the Kronecker delta. Note that $B_\beta$ depends on the braid
\emph{word}, whereas $P_\pi$ only on the braid itself.

\begin{megj}\label{rem:csuszas}
When the braid group relation $\sigma_i\sigma_j=\sigma_j\sigma_i$, $|i-j|>1$
is applied to change $\beta$, the diagram $\gamma_\beta$ only changes by an
isotopy of the plane and the path matrix $B_\beta$ hardly changes at all. In
fact if we don't insist on increasing label indices and re-label the braid as
on the right side of Figure \ref{fig:haromszog}, then $B_\beta$ remains the same.
Therefore such changes in braid words will be largely ignored in the paper.
\end{megj}

\subsection{Multiplicativity}

The path matrix behaves multiplicatively in the following sense: If
two positive braid words $\beta_1$ and $\beta_2$ on $q$ strands are
multiplied as in the braid group (placing $\beta_2$ to the right of
$\beta_1$) to form the braid word $\beta_1*\beta_2$, then
\begin{equation}\label{eq:multip}
B_{\beta_1*\beta_2}=B_{\beta_1}\cdot B_{\beta_2}.
\end{equation}
Note that for this to hold true, $\beta_1$ and $\beta_2$ have to
carry their own individual crossing labels that $\beta_1*\beta_2$
inherits, too. Otherwise, the observation is trivial: we may
group together paths from left endpoint $i$ to right endpoint $j$ in
$\beta_1*\beta_2$ by the position of their crossing over from
$\beta_1$ to $\beta_2$.

\begin{figure}
   \centering
   \includegraphics[width=\linewidth]{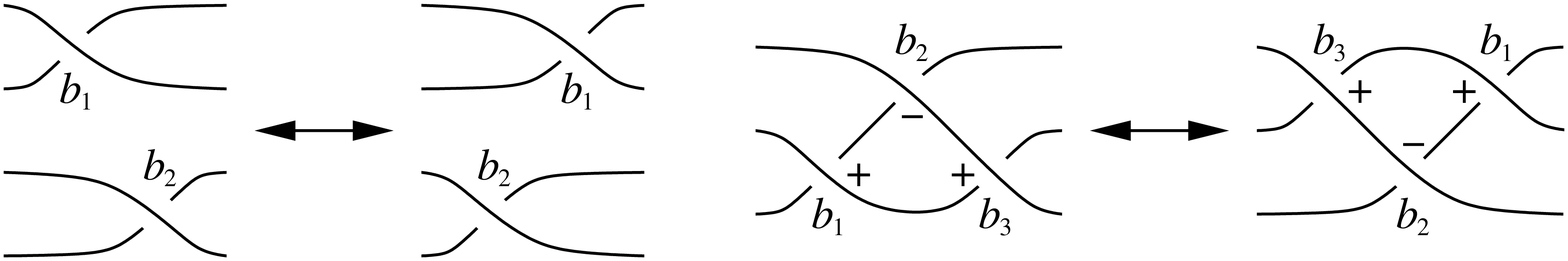}
   \caption{Labels before and after an isotopy and a Reidemeister III move.
   The signs on the right are the so called Reeb signs.}
   \label{fig:haromszog}
\end{figure}

\begin{megj}
Apart from the technicality of having to view $B_{\beta_1}$ and
$B_{\beta_2}$ as polynomials of separate sets of indeterminates, there are
other problems that so far prevented the author from defining a
representation of the positive braid semigroup based on \eqref{eq:multip}.
Namely, when we represent the same positive braid by a different braid
word, the path matrix changes. This can be somewhat controlled by
requiring, as another departure from our convention of increasing label
subscripts, that whenever the braid group relation
$\sigma_i\sigma_{i+1}\sigma_i=\sigma_{i+1}\sigma_i\sigma_{i+1}$ is applied to
change $\beta$, the two
sets of labels are related as on the right side of Figure \ref{fig:haromszog}.
Then the path
matrix changes \begin{equation}\label{eq:haromszog} \text{from
}\begin{bmatrix}b_2&b_3&1\\b_1&1&0\\1&0&0\end{bmatrix}\text{ to
}\begin{bmatrix}b_2+b_3b_1&b_3&1\\b_1&1&0\\1&0&0\end{bmatrix}.\end{equation}
Notice that this is just an application of Chekanov's chain map
\cite{chek} relating the DGA's of the diagrams before and after a
Reidemeister III move (and the same happens if the triangle is part of a
larger braid). Therefore we may hope that the path matrix of a positive
braid $\beta$, with its \emph{entries viewed as elements of the relative
contact homology $H(L_\beta)$}, is independent of the braid word
representing $\beta$. This is indeed the case because the set of equivalent
positive geometric braids (with the endpoints of strands
fixed\footnote{I.e., conjugation is not allowed here; if it was, the space
in question would not be contractible any more, as demonstrated in \cite{en}.})
is contractible, 
thus it is possible
to canonically identify the contact homologies coming from different
diagrams. But because there isn't any known relation between the contact
homologies of $L_{\beta_1}$, $L_{\beta_2}$, and $L_{\beta_1*\beta_2}$,
this doesn't help us. 
\end{megj}

The path matrix of the braid group generator $\sigma_i$, with its
single crossing labeled $b$ is block-diagonal with only two
off-diagonal entries:
\begin{equation}\label{eq:elemi}
B_{\sigma_i}=\left[\begin{array}{rccl}
I_{i-1}&&&\\
&b&1&\\
&1&0&\\
&&&I_{q-i-1}
\end{array}\right].
\end{equation}
By \eqref{eq:multip}, all path matrices are products of such
elementary matrices.

\begin{figure}
   \centering
   \begin{minipage}[c]{.4\textwidth}
   \centering
   \includegraphics[width=\textwidth]{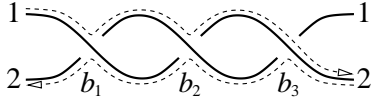}
   \end{minipage}
   \begin{minipage}[c]{.4\textwidth}
   \centering
   \caption{Trefoil braid}
   \label{fig:111}
   \end{minipage}
\end{figure}

\begin{pelda}
Consider the braid $\beta$ shown in Figure \ref{fig:111}. Its path matrix is
\[B_\beta=
\begin{bmatrix}B_{1,1}&B_{1,2}\\B_{2,1}&B_{2,2}\end{bmatrix}=
\begin{bmatrix}b_1&1\\1&0\end{bmatrix}
\begin{bmatrix}b_2&1\\1&0\end{bmatrix}
\begin{bmatrix}b_3&1\\1&0\end{bmatrix}=
\begin{bmatrix}b_1+b_3+b_1b_2b_3&1+b_1b_2\\1+b_2b_3&b_2\end{bmatrix}.\]
(The path contributing $b_1b_2$ to $B_{1,2}$ is shown.)
As $D_1=\{\:\varnothing\:\}$ and
$D_2=\{\:\varnothing,\{\,1\,\}\:\}$, we have
$C_{1,1}=B_{1,1}=b_1+b_3+b_1b_2b_3$ and
$C_{2,2}=B_{2,2}+B_{2,1}B_{1,2}=b_2+(1+b_2b_3)(1+b_1b_2)$. Thus in the DGA of
$\gamma_\beta$, the relations
$\partial a_1=1+b_1+b_3+b_1b_2b_3$ and $\partial
a_2=1+b_2+(1+b_2b_3)(1+b_1b_2)=b_2+b_2b_3+b_1b_2+b_2b_3b_1b_2$ hold.
\end{pelda}


\subsection{Inverse matrix}

The inverse of the elementary matrix $B_{\sigma_i}$ is
\[B^{-1}_{\sigma_i}=\left[\begin{array}{rccl}
I_{i-1}&&&\\
&0&1&\\
&1&b&\\
&&&I_{q-i-1}
\end{array}\right].\]
Therefore, writing
$\beta=\sigma_{i_1}\sigma_{i_2}\cdots\sigma_{i_w}$, from
$B^{-1}_\beta=\left(B_{\sigma_{i_1}}B_{\sigma_{i_2}}\cdots
B_{\sigma_{i_w}}\right)^{-1} = B^{-1}_{\sigma_{i_w}}\cdots
B^{-1}_{\sigma_{i_1}}$ we see that $B^{-1}_\beta$ is also a path
matrix of the same braid word $\beta$, but in a different sense. This
time, the $(i,j)$--entry is a sum of the following products: For
each path composed of parts of the strands of $\beta$ that connects
the \emph{right} endpoint labeled $i$ to the \emph{left} endpoint
labeled $j$ so that it only turns at quadrants facing \emph{down},
take the product of the crossings from right to left that it turns
at. So it's as if we turned $\beta$ upside down by a $180^\circ$ rotation
while keeping the original labels of the crossings
and of the endpoints of the strands.

\begin{figure}
\centering
\includegraphics[width=\linewidth]{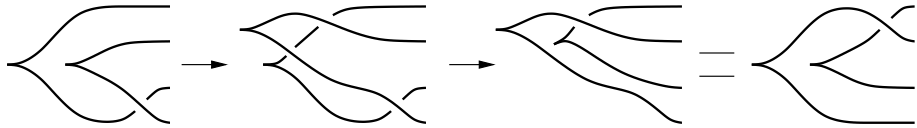}
\caption{First half of a conjugation move (sequence of two Reidemeister II
moves)}
\label{fig:felszalad}
\end{figure}

That operation on the braid word produces a Legendrian isotopic closure
(where by closure we mean adding strands above the braid, as in Figure
\ref{fig:lagpic}). This is seen by a two-step process. First, apply
`half-way' the conjugation move of \cite{en} (as in Figure
\ref{fig:felszalad}) successively to each crossing of $\beta$ from left to
right. This turns $\beta$ upside down, but now the closing strands are
underneath.

\begin{figure}
   \centering
   \includegraphics[width=\linewidth]{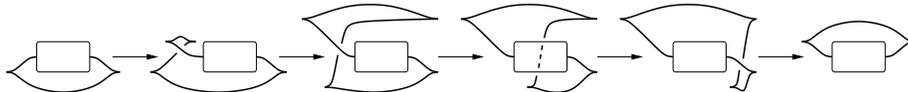}
   \caption{Moving a strand to the other side of a front diagram.}
   \label{fig:felfordul}
\end{figure}

Then, repeat $q$ times the procedure shown in Figure \ref{fig:felfordul},
which we borrow from \cite{tab}. The box may contain any front diagram.
Before the move represented by the third arrow, we make the undercrossing
strand on the left steeper than all slopes that occur inside the box, so
that it slides underneath the entire diagram without a self-tangency
moment. (In $3$--space, increasing the slope results in a huge
$y$--coordinate. Recall that fronts appear on the $xz$--plane, in
particular the $y$--axis points away from the observer. So the motion of
the strand happens far away, way behind any other piece of the knot.)

\begin{pelda}
The inverse of the matrix from the previous example is
\[B^{-1}_\beta=
\begin{bmatrix}0&1\\1&b_3\end{bmatrix}
\begin{bmatrix}0&1\\1&b_2\end{bmatrix}
\begin{bmatrix}0&1\\1&b_1\end{bmatrix}=
\begin{bmatrix}b_2&1+b_2b_1\\1+b_3b_2&b_3+b_1+b_3b_2b_1\end{bmatrix}.\]
In Figure \ref{fig:111}, the path contributing $b_3b_2b_1$ to the $(2,2)$ entry is shown.
\end{pelda}

\subsection{Permutation braids}

As an illustration, we examine the path matrices of permutation
braids, which are positive braids in which every pair
of strands crosses at most once. They are in a one-to-one
correspondence with elements of the symmetric group $S_q$ and they play a crucial role
in Garside's solution \cite{garside} of the word and conjugacy problems in the braid group $B_q$.

It is always possible to represent a braid with a braid word in which the
product $\sigma_i\sigma_{i+1}\sigma_i$ doesn't appear for any $i$. (That
is, all possible triangle moves in which the ``middle strand is pushed
down,'' as in Figure \ref{fig:haromszog} viewed from the right to the
left, have been performed.) Such \emph{reduced braid words} for
permutation braids (up to the relation $\sigma_i\sigma_j=\sigma_j\sigma_i$, 
$|i-j|>1$; see Remark \ref{rem:csuszas}) are unique.


\begin{all}
Let $\pi\in S_q$. The path matrix $B_\pi$ associated to its reduced
permutation braid word is obtained from the permutation matrix
$P_\pi$ as follows. Changes are only made to entries that are above
the $1$ in their column and to the left of the $1$ in their row. At
each such position, a single crossing label appears in $B_\pi$.
\end{all}

In particular, the positions that carry different entries in $P_\pi$
and $B_\pi$ are in a one-to-one correspondence with the inversions
of $\pi$.


\begin{proof} Starting at the left endpoint labeled $i$, our first
``intended destination'' (on the right side of the braid) is $\pi(i)$.
Whenever we turn along a path in the braid, the intended destination
becomes a smaller number because the two strands don't meet again. This
shows that entries in $B_\pi$ that are to the right of the $1$ in their
row are $0$. Traversing the braid from right to left, we see that entries
under the $1$ in their column are $0$, too. Either one of the two
arguments shows that the $1$'s of $P_\pi$ are left unchanged in $B_\pi$.
(This part of the proof is valid for any positive braid word representing
a permutation braid; cf.\ Figure \ref{fig:haromszog} and equation
\eqref{eq:haromszog}.)

We claim that any path in the braid contributing to any $B_{i,j}$ can
contain at most one turn. Assume the opposite: then a strand $s$ crosses
under the strand $t_1$ and then over the strand $t_2$, which are different
and which have to cross each other as well. This contradicts our
assumption that the braid word is reduced, for it is easy to argue that
(in a permutation braid) the triangle $s,t_1,t_2$ that we have just found
must contain an elementary triangle as on the right side of Figure
\ref{fig:haromszog}.

So the
paths we have not yet enumerated are those with exactly one turn.
Because strands cross at most once, these contribute to different
matrix entries. Finally, if $(i,j)$ is a position as described in
the Proposition, then $\pi(i)>j$ and $\pi^{-1}(j)>i$. This means
that the strand starting at $i$ has to meet the strand ending at
$j$, so that the label of that crossing becomes $B_{i,j}$.
\end{proof}


\begin{figure}
   \centering
   \includegraphics[width=.9\linewidth]{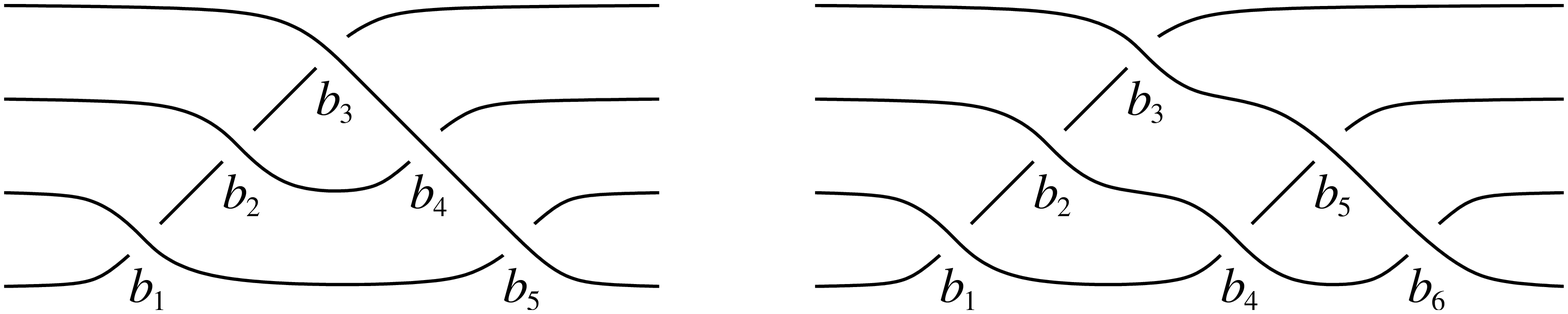}
   \caption{Permutation braids of $(14)$ and of $(14)(23)$
   (the latter also known as the Garside braid $\Delta_4$).}
   \label{fig:permbraids}
\end{figure}

\begin{pelda}
The transposition $(14)$ of $S_4$ is represented by the reduced braid word shown
in Figure \ref{fig:permbraids}. It contains $5$ inversions, corresponding to the
$5$ crossings of the braid. Its path matrix is
$B_{(14)}=\begin{bmatrix}b_3&b_4&b_5&1\\b_2&1&0&0\\b_1&0&1&0\\1&0&0&0
\end{bmatrix}$.
The path matrix of the Garside braid (half-twist) $\Delta_4$, also
shown in Figure \ref{fig:permbraids}, is
$\begin{bmatrix}b_3&b_5&b_6&1\\b_2&b_4&1&0\\b_1&1&0&0\\1&0&0&0\end{bmatrix}$.
\end{pelda}

The latter pattern obviously generalizes to $\Delta_n$ for any $n$.

\subsection{Row reduction}

There is yet another way to factorize the path matrix. Let
$\tau_i\in S_q$ denote the underlying permutation (transposition) of
the elementary braid $\sigma_i\in B_q$.

\begin{lemma}\label{lem:atmegy}
Let $\lambda\in S_q$ be an arbitrary permutation. Then for all $i$,
\begin{equation}\label{eq:atmegy}\left[\begin{array}{c}
\text{matrix}\\
\text{of }\lambda
\end{array}\right]\cdot
\left[\begin{array}{rccl}
I_{i-1}&&&\\
&b&1&\\
&1&0&\\
&&&I_{q-i-1}
\end{array}\right]=
\left[\begin{array}{rccl}1&&&\\&1&&\\&b&\ddots&\\&&&1\end{array}\right]\cdot
\left[\begin{array}{c}
\text{matrix}\\
\text{of }\tau_i\circ\lambda\end{array}\right],\end{equation}
where in the first
term of the right hand side, the single non-zero off-diagonal entry
$b$ appears in the position $\lambda^{-1}(i),\lambda^{-1}(i+1)$.
\end{lemma}

\begin{figure}
   \centering
   \includegraphics[width=4in]{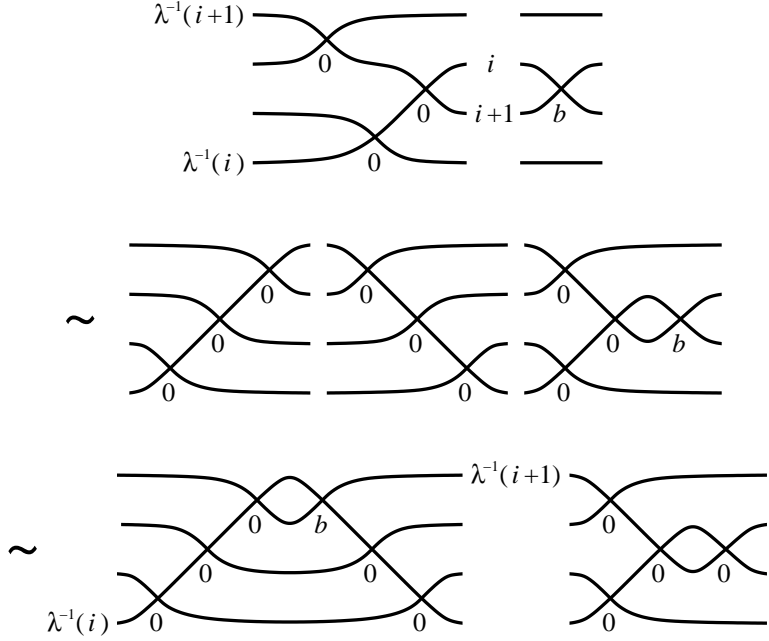}
   \caption{Braids so decorated that they have the same path matrix}
   \label{fig:bizony}
\end{figure}

\begin{proof}
The essence of the proof is in Figure
\ref{fig:bizony}. It will be crucial that the path matrix depends on how
the braid is decorated with labels. On the other hand, for the purposes of
the argument, over- and
undercrossing information in the braids is irrelevant. In fact, although
we will not change our terminology, we will actually think of them (in
particular, when we take an inverse) as words written in the generators
$\tau_1,\ldots,\tau_{q-1}$ of $S_q$.

Take the permutation braid for $\lambda$ (or choose any other positive braid word
with this underlying permutation) and label its crossings with zeros. (In
Figure \ref{fig:bizony} we used $\lambda=(1342)\in S_4$ as an example.)
Add a single generator $\sigma_i$, with its crossing labeled $b$ to it
(Figure \ref{fig:bizony} shows $i=2$). The left hand side of
\eqref{eq:atmegy} is the path matrix of this braid $\beta$.

Next, choose any positive braid word $\mu$ in which the strands with right
endpoints $\lambda^{-1}(i)$, $\lambda^{-1}(i+1)$ cross (say exactly once)
and form the product $\mu^{-1}*\mu*\beta$. Label the crossings of
$\mu^{-1}$ and $\mu$ with zeros, as in the middle of Figure
\ref{fig:bizony}. This way, the path matrix does not change.

Now, it does not matter for the path matrix where exactly the single
non-zero label $b$ appears in the braid as long as that crossing
establishes a path between the same two endpoints. In other words, we may
move the label from the first (from the right) to the third, fifth etc.\
crossing of the same two strands. By construction, one of those crossings
is either in $\mu$ (if $\lambda^{-1}(i)>\lambda^{-1}(i+1)$, as is the case
in Figure \ref{fig:bizony}) or in $\mu^{-1}$, and we move the label there
(bottom of Figure \ref{fig:bizony}). When we read off the path matrix from
this form, we obtain the right hand side of \eqref{eq:atmegy}: The path
matrix of $\mu^{-1}*\mu$ is $I_q$ except for the single $b$ that
establishes a path from $\lambda^{-1}(i)$ to $\lambda^{-1}(i+1)$, and the
path matrix of $\beta$, now labeled with only zeros, is
$P_{\tau_i\circ\lambda}$.
\end{proof}

Next, for the positive braid word
$\beta=\sigma_{i_1}\sigma_{i_2}\cdots\sigma_{i_w}$ with crossings labeled
$b_1,b_2,\ldots,b_w$, we'll introduce a sequence of elementary matrices.
The underlying permutation is $\pi=\tau_{i_w}\ldots\tau_{i_2}\tau_{i_1}$.
Let us denote the ``permutation up to the $k$'th crossing'' by
$\pi_k=\tau_{i_k}\ldots\tau_{i_1}$, so that $\pi_0=\id$ and
$\pi_w=\pi$. Let $A_k$ be the $q\times q$ identity matrix with a single
non-zero off-diagonal entry of $b_k$ added in the position
$\pi_{k-1}^{-1}(i_k),\pi_{k-1}^{-1}(i_k+1)$. Note that because we work over
$\Z_2$, $A_k^2=I_q$ for all $k$.

\begin{all}
For the positive braid word
$\beta=\sigma_{i_1}\sigma_{i_2}\cdots\sigma_{i_w}$ with underlying
permutation $\pi$, we have $B_\beta=A_1A_2\ldots A_wP_\pi$, where
$P_\pi$ is the permutation matrix.
\end{all}

\begin{proof} If $\lambda=\pi_{k-1}$, $i=i_k$, and $b=b_k$, then equation
\eqref{eq:atmegy} reads $P_{\pi_{k-1}} B_{\sigma_{i_k}} = A_k
P_{\pi_k}$. Starting from $B_\beta=(P_{\pi_0}
B_{\sigma_{i_1}})B_{\sigma_{i_2}}\ldots B_{\sigma_{i_w}}$, we apply Lemma
\ref{lem:atmegy} $w$ times.
\end{proof}

Read in another way, this result shows that $B_\beta$ reduces to $P_\pi$ by
applying a particular sequence of elementary row operations:
$A_w\ldots A_1B_\beta=P_\pi$. This works in the non-commutative sense.

\begin{pelda} For the braid $\beta$ of Figure \ref{fig:111}, we have
\[B_\beta=\begin{bmatrix}1&b_1\\0&1\end{bmatrix}
\begin{bmatrix}1&0\\b_2&1\end{bmatrix}
\begin{bmatrix}1&b_3\\0&1\end{bmatrix}
\begin{bmatrix}0&1\\1&0\end{bmatrix},\]
that is
\[\begin{bmatrix}1&b_3\\0&1\end{bmatrix}
\begin{bmatrix}1&0\\b_2&1\end{bmatrix}
\begin{bmatrix}1&b_1\\0&1\end{bmatrix}
\begin{bmatrix}b_1+b_3+b_1b_2b_3&1+b_1b_2\\1+b_2b_3&b_2\end{bmatrix}
=\begin{bmatrix}0&1\\1&0\end{bmatrix}.\]
\end{pelda}

\section{Algebraic results}

In this section, we treat (re-define, if you like) the symbols
$B_{i,j}$ as independent variables. Instead of $\Z_2$--coefficients,
we will work in the free non-commutative unital ring generated by
these symbols (where $1\le i,j\le q$) \emph{over $\Z$}.
After the first set
of statements, we will abelianize so that we can consider
determinants.

Note that the $C_{i,j}$ (equation \eqref{eq:Cij}) are polynomials in the
$B_{i,j}$. To state our results, we will need a similar family of
polynomials whose definition is based on the notion of admissible sequence
(Definition \ref{def:sorozat}).

\begin{Def} For any $1\le i,j\le q$, let
\begin{equation*}
M_{i,j}=\sum_{\{\,i_1,\ldots,i_c\,\}\in D_{\min\{i,j\}}}
B_{i,i_1}B_{i_1,i_2}B_{i_2,i_3}\ldots B_{i_{c-1},i_c}B_{i_c,j}.
\end{equation*}
\end{Def}

Note that $M_{1,j}=C_{1,j}=B_{1,j}$, $M_{i,1}=B_{i,1}$, $M_{n,n}=C_{n,n}$, and
$M_{i-1,i}=C_{i-1,i}$,
whenever these expressions are defined.

\begin{lemma}
\begin{multline*}\left[\begin{array}{ccccc}
1&C_{1,2}&C_{1,3}&\cdots&C_{1,q}\\
&1&C_{2,3}&\cdots&C_{2,q}\\
&&1&\cdots&C_{3,q}\\
&&&\ddots&\vdots\\
&&&&1
\end{array}\right]\cdot
\left[\begin{array}{ccccc}
1&-M_{1,2}&-M_{1,3}&\cdots&-M_{1,q}\\
&1&-M_{2,3}&\cdots&-M_{2,q}\\
&&1&\cdots&-M_{3,q}\\
&&&\ddots&\vdots\\
&&&&1
\end{array}\right]\\
= \left[\begin{array}{ccccc}
1&-M_{1,2}&-M_{1,3}&\cdots&-M_{1,q}\\
&1&-M_{2,3}&\cdots&-M_{2,q}\\
&&1&\cdots&-M_{3,q}\\
&&&\ddots&\vdots\\
&&&&1
\end{array}\right]\cdot
\left[\begin{array}{ccccc}
1&C_{1,2}&C_{1,3}&\cdots&C_{1,q}\\
&1&C_{2,3}&\cdots&C_{2,q}\\
&&1&\cdots&C_{3,q}\\
&&&\ddots&\vdots\\
&&&&1
\end{array}\right]
=I_q,\end{multline*} and a similar statement can be formulated for lower
triangular matrices.
\end{lemma}

Note that the two claims don't imply each other because we work
over a non-commutative ring.

\begin{proof}
We need that for all $1\le i<j\le q$,
\[-M_{i,j}-C_{i,i+1}M_{i+1,j}-C_{i,i+2}M_{i+2,j}-\ldots
-C_{i,j-1}M_{j-1,j}+C_{i,j}=0\]
and that
\[C_{i,j}-M_{i,i+1}C_{i+1,j}-M_{i,i+2}C_{i+2,j}-\ldots
-M_{i,j-1}C_{j-1,j}-M_{i,j}=0.\]
We may view both of these equalities as identities for $C_{i,j}$.
The first one groups the terms of $C_{i,j}$ according to the highest
element of the admissible sequence. The second groups them according
to the first element which is greater than $i$. The lower triangular
version is analogous.
\end{proof}

\begin{lemma}\label{lem:MMB}
For all $1\le n\le q$,
\begin{multline}\label{eq:MMB}
\left[\begin{array}{ccccc}
-1&&&&\\
M_{2,1}&-1&&&\\
M_{3,1}&M_{3,2}&-1&&\\
\vdots&\vdots&\vdots&\ddots&\\
M_{n,1}&M_{n,2}&M_{n,3}&\cdots&-1
\end{array}\right]\cdot
\left[\begin{array}{ccccc}
1&-M_{1,2}&-M_{1,3}&\cdots&-M_{1,n}\\
&1&-M_{2,3}&\cdots&-M_{2,n}\\
&&1&\cdots&-M_{3,n}\\
&&&\ddots&\vdots\\
&&&&1
\end{array}\right]\\
=\left[\begin{array}{cccccc}
-1&B_{1,2}&\cdots&B_{1,i}&\cdots&B_{1,n}\\
B_{2,1}&-1-B_{2,1}B_{1,2}&\cdots&B_{2,i}&\cdots&B_{2,n}\\
\vdots&\vdots&\ddots&\vdots&&\vdots\\
B_{i,1}&B_{i,2}&\cdots&B_{i,i}-C_{i,i}-1&\cdots&B_{i,n}\\
\vdots&\vdots&&\vdots&\ddots&\vdots\\
B_{n,1}&B_{n,2}&\cdots&B_{n,i}&\cdots&B_{n,n}-C_{n,n}-1
\end{array}\right].
\end{multline}
\end{lemma}

\begin{proof}
For entries above the diagonal ($i<j$), the claim is that
\[B_{i,j}=-M_{i,1}M_{1,j}-M_{i,2}M_{2,j}-\ldots-M_{i,i-1}M_{i-1,j}+M_{i,j}.\]
Viewing this as an identity for $M_{i,j}$, we see that it holds
because terms are grouped with respect to the highest element in the
admissible sequence. The reasoning is the same for positions below
the diagonal. For the diagonal entries, we need to show that
\[B_{i,i}-C_{i,i}-1=-M_{i,1}M_{1,i}-M_{i,2}M_{2,i}-\ldots
-M_{i,i-1}M_{i-1,i}-1.\]
Isolating $C_{i,i}$ this time, we again see a separation of its
terms according to the highest element of the admissible sequence.
\end{proof}

For the rest of the section, we will work in the \emph{commutative}
polynomial ring generated over $\Z$ by the $B_{i,j}$, so that we can
talk about determinants.

\begin{tetel}\label{thm:ideal}
The ideal $I'$ generated by the polynomials
\[1+C_{1,1},\quad 1+C_{2,2},\quad\ldots,\quad 1+C_{q,q}\]
agrees with the ideal $I$
generated by the polynomials
\[L_1=B_{1,1}+1,\enskip L_2=\begin{vmatrix}B_{1,1}&B_{1,2}\\B_{2,1}&B_{2,2}
\end{vmatrix}-1,\enskip \ldots,\enskip
L_q=\begin{vmatrix}B_{1,1}&\cdots&B_{1,q}\\
\vdots&\ddots&\vdots\\
B_{q,1}&\cdots&B_{q,q}\end{vmatrix}-(-1)^q.\]
\end{tetel}

\begin{proof}
Let $n\le q$ and take determinants of both sides of equation
\eqref{eq:MMB}: $(-1)^n$, on the left hand side, agrees with
$\begin{vmatrix}B_{1,1}&\cdots&B_{1,n}\\
\vdots&\ddots&\vdots\\
B_{n,1}&\cdots&B_{n,n}\end{vmatrix}$ plus an element of $I'$ on the
right hand side. Thus, $L_n\in I'$ for all $n$.

The proof of the other containment relation is also based on
equation \eqref{eq:MMB} and goes by induction on $n$. Note that
$1+C_{1,1}=L_1$ and assume that
$1+C_{1,1},1+C_{2,2},\ldots,1+C_{n-1,n-1}$ are all in $I$ (actually,
they are in the ideal generated by $L_1,L_2,\ldots,L_{n-1}$).
Re-writing the determinant of the matrix on the right hand side of
\eqref{eq:MMB}, we find that
\begin{multline*}\hspace{-8pt}(-1)^n=\left|\begin{array}{cccccc}
-1&B_{1,2}&\cdots&B_{1,n-1}&B_{1,n}\\
B_{2,1}&-1-B_{2,1}B_{1,2}&\cdots&B_{2,n-1}&B_{2,n}\\
\vdots&\vdots&\ddots&\vdots&\vdots\\
B_{n-1,1}&B_{n-1,2}&\cdots&B_{n-1,n-1}-C_{n-1,n-1}-1&B_{n-1,n}\\
B_{n,1}&B_{n,2}&\cdots&B_{n,n-1}&B_{n,n}
\end{array}\right|\\
-\left|\begin{array}{cccccc}
-1&B_{1,2}&\cdots&B_{1,n-1}&B_{1,n}\\
B_{2,1}&-1-B_{2,1}B_{1,2}&\cdots&B_{2,n-1}&B_{2,n}\\
\vdots&\vdots&\ddots&\vdots&\vdots\\
B_{n-1,1}&B_{n-1,2}&\cdots&B_{n-1,n-1}-C_{n-1,n-1}-1&B_{n-1,n}\\
0&0&\cdots&0&1+C_{n,n}
\end{array}\right|.
\end{multline*}
Notice that the second determinant is $(-1)^{n-1}(1+C_{n,n})$ (by
Lemma \ref{lem:MMB}), while the first is
$\begin{vmatrix}B_{1,1}&\cdots&B_{1,n}\\
\vdots&\ddots&\vdots\\
B_{n,1}&\cdots&B_{n,n}\end{vmatrix}$ plus an element of $I'$, but
the latter, by the inductive hypothesis, is also in $I$. Isolating
$1+C_{n,n}$, we are done.
\end{proof}


So we see that the ideal $I$ defined in terms of the upper left
corner subdeterminants of the general determinant is also generated
by the polynomials $1+C_{n,n}$, which arise from contact homology
(counting holomorphic discs).
In fact much more is true: the $1+C_{n,n}$ form the reduced
Gr\"obner basis for $I$. Of course this can only be true for certain
term orders that we'll describe now.

In the (commutative) polynomial ring $\Z[B_{i,j}]$, take any order $\prec$
of the indeterminates where any diagonal entry $B_{i,i}$ is larger than
any off-diagonal one. Extend this order to the monomials
lexicographically. (But not degree lexicographically! For example,
$B_{2,2}\succ B_{2,1}B_{1,2}$.) This is a multiplicative term order.

\begin{tetel}\label{thm:grob} The polynomials $1+C_{n,n}$, $n=1,\ldots,q$ (defined in
equation \eqref{eq:Cii}), form the reduced Gr\"obner basis for the ideal
\[I=\left\langle
B_{1,1}+1,\quad\begin{vmatrix}B_{1,1}&B_{1,2}\\B_{2,1}&B_{2,2}\end{vmatrix}-1,
\quad\ldots,\quad \begin{vmatrix}B_{1,1}&\cdots&B_{1,q}\\
\vdots&\ddots&\vdots\\
B_{q,1}&\cdots&B_{q,q}\end{vmatrix}-(-1)^q\right\rangle\] under any of the
term orders $\prec$ described above. \end{tetel}

\begin{proof}
This is obvious from the definitions (see for example \cite{bernd}),
after noting that the initial term of $1+C_{n,n}$ is $B_{n,n}$ and
that by the definition of an admissible sequence, no other term in
$1+C_{n,n}$ contains any $B_{i,i}$. (The initial ideal of $I$ is that
generated by the $B_{n,n}$.)
\end{proof}


\section{Augmentations}

\begin{Def}\label{def:aug} Let $\gamma$ be a Lagrangian diagram of a
Legendrian link $L$. If $L$ has more than one components, we assume that
an admissible grading of the DGA of $\gamma$ has been chosen, too. An
\emph{augmentation} is a subset $X$ of the crossings (the \emph{augmented
crossings}) of $\gamma$ with the following properties.
\begin{itemize}
\item The index of each element of $X$ is $0$.
\item For each generator
$a$ of index $1$, the number of admissible discs with positive corner $a$
and all negative corners in $X$ is even.
\end{itemize}
\end{Def}

Here, an admissible disc is the central object of Chekanov--Eliashberg
theory: These discs determine the differential $\partial$ of the DGA
$\mathscr A$, and thus contact homology $H(L)$. Unlike most of the
literature, we expand the notion of augmentation here (in the
multi-component case) by allowing `mixed' crossings between different
components to be augmented, as long as they have index $0$ in the one
grading we have chosen. Such sets of crossings would typically not be
augmentations for other admissible gradings because it's exactly the index
of a mixed crossing that is ambiguous. Our motivation is that
$\gamma_\beta$, even if it is of multiple components, has the natural
admissible grading introduced in section \ref{sec:prelim}.

The evaluation homomorphism (which is defined on the link DGA, and which
is also called an augmentation) $\varepsilon_X\colon\mathscr A\to\Z_2$
that sends elements of $X$ to $1$ and other generators to $0$, gives rise
to an algebra homomorphism $(\varepsilon_X)_*\colon H(L)\to\Z_2$. In fact,
the second requirement of Definition \ref{def:aug} is just an elementary
way of saying that $\varepsilon_X$ vanishes on
$\partial(a)$ for each generator $a$ of index $1$, while for other indices
this is already automatic by the first point and the fact that $\partial$
lowers the index by $1$.

\begin{megj} As a preview of a forthcoming paper, let us mention that
augmentations do define a Legendrian isotopy invariant in the following
sense: the set of all induced maps $(\varepsilon_X)_*\colon H(L)\to\Z_2$
depends only on $L$. (The correspondence between augmentations of
different diagrams of $L$ is established using pull-backs by the
isomorphisms constructed in Chekanov's proof of the invariance of $H(L)$.)
The number of augmentations in the sense of Definition \ref{def:aug} may however
change by a factor of $2$ when a Reidemeister II move or its inverse,
involving crossings of index $0$ and $-1$, is performed. \end{megj}



In practice, finding an augmentation means solving a system of polynomial equations
(one equation provided by each index $1$ crossing) over $\Z_2$. In this sense,
augmentations form a variety. In this section we prove a few statements about the
variety associated to $\gamma_\beta$.

The main result is the following
theorem, which allows for an enumeration of all
augmentations of $\gamma_\beta$. The author is greatly indebted to Supap Kirtsaeng,
who wrote a computer program based on this criterion. It may first seem ineffective
to check all subsets of the crossings of $\beta$, but it turns out that a significant
portion of them are augmentations (see section \ref{sec:ex}).

Let $Y$ be a subset of the crossings of $\beta$. Let $\varepsilon_Y\colon\mathscr A\to\Z_2$
be the evaluation homomorphism that sends elements of $Y$ to $1$ and other generators to $0$.
In particular, we may talk of the $0$-$1$--matrix $\varepsilon_Y(B_\beta)$. (This could also
have been denoted by $B_\beta(\chi_Y)$, where the $0$-$1$--sequence $\chi_Y$ is the
characteristic function of $Y$.)

\begin{tetel}\label{thm:gauss}
Let $Y$ be a subset of the crossings of the positive braid word $\beta$. $Y$ is an
augmentation of $\gamma_\beta$ if and only if the
$0$-$1$--matrix $\varepsilon_Y(B_\beta)$ is such that every upper left
corner square submatrix of it has determinant $1$.
\end{tetel}

It is then a classical theorem of linear algebra that the condition
on $\varepsilon_Y(B_\beta)$ is equivalent to the requirement that it
possess an $LU$--decomposition and also to the requirement that
Gaussian elimination can be completed on it without permuting rows.

\begin{proof} In our admissible grading, each crossing of $\beta$ has
index $0$. Therefore $Y$ is an augmentation if and only if $\varepsilon_Y$
vanishes on $\partial(a)$ for each index $1$ DGA generator $a$.  This in
turn is clearly equivalent to saying that $\varepsilon_Y$ vanishes on the
two-sided ideal generated by these polynomials. In fact because
$\varepsilon_Y$ maps to a commutative ring ($\Z_2$), we may abelianize
$\mathscr A$ and say that the condition for $Y$ to be an augmentation is
that $\varepsilon_Y$ vanishes on the ideal generated by the expressions
$\partial(a_1),\ldots,\partial(a_q)$, which are now viewed as honest
polynomials in the commuting indeterminates $b_1,\ldots,b_w$.

In \cite{en}, section 6, we computed these polynomials and found that they
really were polynomials of the polynomials $B_{i,j}$, as stated in
equation \eqref{eq:perem}. Now by (the modulo $2$ reduction of) Theorem
\ref{thm:ideal}, the ideal generated by the $\partial(a_n)$ is also
generated by the polynomials
\[B_{1,1}+1,\quad\begin{vmatrix}B_{1,1}&B_{1,2}\\B_{2,1}&B_{2,2}\end{vmatrix}
+1,\quad\ldots,\quad
\begin{vmatrix}B_{1,1}&\cdots&B_{1,q}\\ \vdots&\ddots&\vdots\\
B_{q,1}&\cdots&B_{q,q}\end{vmatrix}+1,\] which implies the Theorem
directly.
\end{proof}

\begin{megj}\label{rem:utso} Notice that for a path matrix $B_\beta$, a
quick look at \eqref{eq:elemi} with formula \eqref{eq:multip} implies that
we always have $\det(B_\beta)=1$. Therefore the condition on the $q\times
q$ subdeterminant is vacuous: if a subset of the crossings of $\beta$
``works as an augmentation'' for $a_1,\ldots,a_{q-1}$, then it
automatically works for $a_q$ as well. \end{megj}

\begin{figure}
   \centering
   \includegraphics[width=4in]{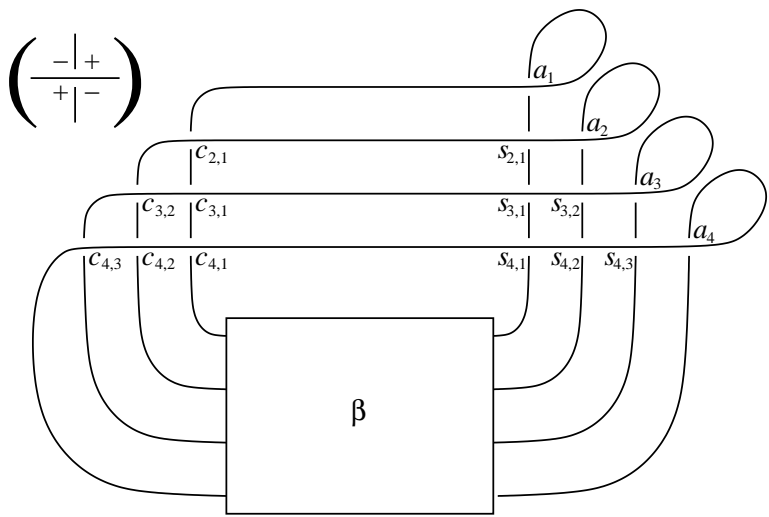}
   \caption{Another Lagrangian diagram of $L_\beta$}
   \label{fig:plat}
\end{figure}

Let us give a geometric explanation of the appearance of
$LU$--de\-com\-po\-si\-tions. Figure \ref{fig:plat} shows another Lagrangian
diagram of $L_\beta$ that is obtained from the front diagram $f_\beta$ by
pushing all the right cusps to the extreme right and then applying resolution.
This has the advantage that all admissible discs are embedded. Label the
$q(q-1)$ new crossings as in Figure \ref{fig:plat}. Our preferred grading is
extended to the new crossings by assigning $0$ to the $c_{i,j}$ and $1$ to the
$s_{i,j}$. This implies $\partial(c_{i,j})=0$, while the index $1$ generators
are mapped as follows: \[\partial(a_n) = 1 + c_{n,1}B_{1,n} +\ldots
+c_{n,n-1}B_{n-1,n} + B_{n,n}\] and \[\partial(s_{i,j}) = c_{i,1}B_{1,j}
+\ldots+ c_{i,i-1}B_{i-1,j} + B_{i,j}.\] Setting the latter $q+{q\choose 2}$
expressions equal to $0$ is equivalent to saying that the matrix product
\[\left[\begin{array}{ccccc}
1&0&0&\cdots&0\\
c_{2,1}&1&0&\cdots&0\\
c_{3,1}&c_{3,2}&1&\cdots&0\\
\vdots&\vdots&\vdots&\ddots&\vdots\\
c_{q,1}&c_{q,2}&c_{q,3}&\cdots&1
\end{array}\right]
\left[\begin{array}{ccccc}
B_{1,1}&B_{1,2}&B_{1,3}&\cdots&B_{1,q}\\
B_{2,1}&B_{2,2}&B_{2,3}&\cdots&B_{2,q}\\
B_{3,1}&B_{3,2}&B_{3,3}&\cdots&B_{3,q}\\
\vdots&\vdots&\vdots&\ddots&\vdots\\
B_{q,1}&B_{q,2}&B_{q,3}&\cdots&B_{q,q}
\end{array}\right]\]
is unit upper triangular. Thus an augmentation evaluates $B_\beta$ to an
$LU$--decomposable $0$-$1$--matrix and the converse is not hard to prove either.

\section{Rulings}\label{sec:rul}

\begin{figure}
   \centering
   \includegraphics[width=\linewidth]{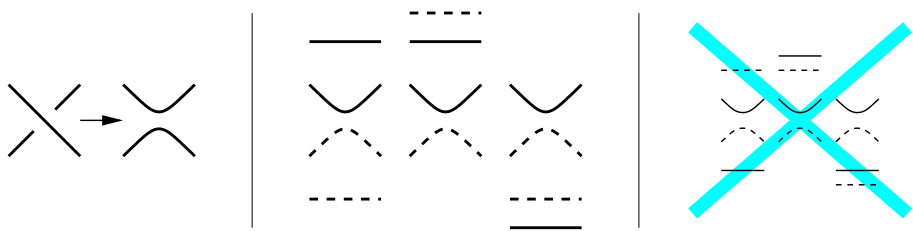}
   \caption{Allowed and disallowed configurations for switches of rulings}
   \label{fig:switch}
\end{figure}

\begin{Def}\label{def:rul} An \emph{ungraded ruling} is a partial splicing
of a front diagram where certain crossings, called \emph{switches}, are
replaced by a pair of arcs as in Figure \ref{fig:switch} so that the
diagram becomes a (not necessarily disjoint) union of standard unknot diagrams, called \emph{eyes}.
(An eye is a pair of arcs connecting the same two cusps that contain no
other cusps and that otherwise do not meet, not even at switches.) It is
assumed that in the vertical ($x=\text{const.}$) slice of the diagram
through each switch, the two eyes that meet at the switch follow one of
the three configurations in the middle of Figure \ref{fig:switch}.

Let us denote the set of all ungraded rulings of a front diagram $f$ of a Legendrian
link by $\Gamma_1(f)$. We get $2$--graded rulings, forming the set $\Gamma_2(f)$, if
we require that the index of each switch be
even. $\Z$--graded rulings (set $\Gamma_0(f)$) are those where each switch has index $0$.
\end{Def}

$\Gamma_1$ is of course grading-independent. For multi-component oriented link diagrams, 
$\Gamma_2$ doesn't depend on the chosen grading, but $\Gamma_0$ might. 

Rulings can also be classified by the value
\[\theta=\text{number of eyes}-\text{number of switches}.\] The counts of
ungraded, $2$--graded, and $\Z$--graded rulings with a given $\theta$ are
all Legendrian isotopy invariants\footnote{In the $\Z$--graded case, we may have to assume that the Legendrian has a single component.} \cite{chp,fuchs}. (In particular, the sizes of the
sets $\Gamma_i(f)$, $i=0,1,2$, don't depend on $f$, only on the Legendrian isotopy class.)
We may arrange these numbers as coefficients in the ruling polynomials\footnote{These are
honest polynomials for knots, but for multi-component links, they may contain negative
powers of $z$. It may seem unnatural first to write them the way we do, but there are two
good reasons to do so: One is Rutherford's pair of theorems below, and the other is that
rulings can also be thought of as surfaces, in which case $\theta$ becomes their Euler
characteristic and (in the one-component and $2$--graded case) $1-\theta$ is twice their genus.}
\[R_i(z)=\sum_{\rho\in\Gamma_i}z^{1-\theta(\rho)}.\]

Fuchs
notes that the existence of a $2$--graded ruling
implies $r=0$. Let us add that if we treat the eyes as discs and
join them by twisted bands at the switches, then a $2$--graded
ruling becomes an orientable surface. The number $\theta$ is its
Euler characteristic and thus $\theta+\mu$, where $\mu$ is the
number of the components of the Legendrian, is even. In particular,
$\theta$ is odd for any $2$--graded ruling of a Legendrian knot.

There is a marked difference between $\Z$--graded rulings and the two less restrictive
cases. $R_1$ and $R_2$ only depend on the smooth type of the Legendrian and its
Thurston--Bennequin number. In fact, Rutherford \cite{rulpoly} proved that for any link,
$R_1(z)$ is the coefficient of $a^{-tb-1}$ in the Dubrovnik version of the Kauffman
polynomial, and that $R_2(z)$ is the coefficient of $v^{tb+1}$ in the Homfly polynomial.
On the other hand, $R_0$ is more sensitive: Chekanov \cite{chek2} constructed two
Legendrian knots of type $5_2$, both with $tb=1$ and $r=0$, so that one has $R_0(z)=1+z^2$
and the other has $R_0(z)=1$.

Because $f_\beta$ only contains crossings of index $0$, any ungraded
ruling is automatically $2$--graded and $\Z$--graded in this case. Thus we may talk about
a single ruling polynomial. By Rutherford's theorems, this implies that the coefficients
of the terms with minimum $v$--degree in the Homfly and Kauffman
polynomials (for the latter, replace $a$ with $v^{-1}$ in its Dubrovnik version) of a
braid-positive link agree. In fact, using Tanaka's results \cite{tanaka}, the same can be said about arbitrary positive links. (See \cite{meginten} for more.) This, without any reference to Legendrians yet with essentially the same proof, has been first observed by Yokota \cite{yok}.

\begin{pelda} The positive trefoil knot that is the closure of the braid
in Figure \ref{fig:111} has one ruling with $\theta=-1$ and two with
$\theta=1$, shown in Figure \ref{fig:rulings}. The numbers $1$ and $2$ (i.e., the ruling polynomial $R(z)=2+z^2$) 
appear as the leftmost coefficients in the Homfly polynomial
\[\begin{array}{ccc}z^2v^2&&\\&&\\2v^2&&-v^4\end{array}\] and also in the Kauffman polynomial
\[\begin{array}{cccc}z^2v^2&&-z^2v^4&\\&-zv^3&&+zv^5\\2v^2&&-v^4.\end{array}\]

\begin{figure}
   \centering
   \includegraphics[width=\linewidth]{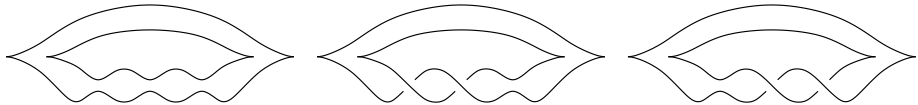}
   \caption{The Seifert ruling and the other two rulings of the positive trefoil}
   \label{fig:rulings}
\end{figure}
\end{pelda}

The diagram $f_\beta$ admits many rulings. The one that is easiest to see
is what we will call the \emph{Seifert ruling}, in which the set of
switches agrees with the set of crossings in $\beta$. This is the only
ruling with the minimal value $\theta=q-w$. Another ruling, that one with the maximum
value $\theta=\mu$, will be constructed in Theorem \ref{thm:szimultan}.

The second lowest possible value of $\theta$ for a ruling of $f_\beta$ is
$q-w+2=-tb(L_\beta)+2$. It is easy to see that in such a ruling, the two crossings
of $\beta$ that are not switches have to be `on the same level' (represented by the
same braid group generator) without any other crossing between them on that level, and
also that any such arrangement works. Thus, assuming that each generator occurs in the
braid word $\beta$ (i.e., that $f_\beta$ is connected), the number of such rulings is
$w-(q-1)=tb(L_\beta)+1$. In all of the examples known to the author, the next value of
$\theta$, that is $\theta=q-w+4=-tb+4$, is realized by exactly
${{w-q}\choose{2}}={{tb}\choose{2}}$ rulings (but I don't know how to prove this).
At $\theta=-tb+6$ and higher, dependence on the braid occurs (see section \ref{sec:ex}).

It would be very interesting to have a test, similar to Theorem \ref{thm:gauss}, that
decides from the path matrix whether a given crossing set of $f_\beta$ is a ruling.

From work of Fuchs, Ishkhanov \cite{fuchs, masikirany}, and Sabloff \cite{josh}, we
know that $\Z$--graded
rulings for a Legendrian exist if and only if augmentations do.
Ng and Sabloff also worked out a
surjective correspondence \cite{manytoone} that assigns a $\Z$--graded ruling to each
augmentation. In that correspondence, the size of the preimage of each $\Z$--graded
ruling $\rho$ of the front diagram $f$ is the number $2^{(\theta(\rho)+\chi^*(f))/2}$, where
\begin{eqnarray*}
\chi^*(f)&\hspace{-8pt}
=&\hspace{-8pt}- \left(\sum_{\text{crossings } a \text{ of } f \text{ with } |a|<0}(-1)^{|a|}
\right)
+ \left(\sum_{\text{crossings } a \text{ of } f \text{ with } |a|\ge 0}(-1)^{|a|}\right)\\
&&\hspace{-8pt}-\hspace{5pt}\text{number of right cusps}.
\end{eqnarray*}
In particular, the number of augmentations belonging to $\rho$ depends on $\theta(\rho)$ and
the diagram only. (Note that $\chi^*$ has the same parity as $tb$, and because $r=0$ is
even, it also has the same parity as $\mu$.) Thus the total number of augmentations is
\begin{equation}\label{eq:aug}
R_0(z)\cdot z^{-1-\chi^*}\bigg|_{z=2^{-1/2}}.
\end{equation}

For the diagram $f_\beta$, which is without negatively graded crossings, we have
$\chi^*(f_\beta)=tb(L_\beta)=w-q$. Thus among the rulings of
$f_\beta$, the zeroth power of $2$ corresponds only to the Seifert ruling.
Therefore the number of augmentations of $f_\beta$ is odd.

The next theorem may further illuminate the relationship between augmentations and rulings.

\begin{tetel}\label{thm:szimultan}
For any positive braid word $\beta$, there exists a subset of
its crossings which is (the set of switches in) a ruling of $f_\beta$ and
an augmentation of $\gamma_\beta$ at the same time.
\end{tetel}

The set we will construct is not, however, fixed by Ng and Sabloff's
many-to-one correspondence.

\begin{proof} The set $X$ is constructed as follows: In $\beta$, the
strands starting at the left endpoint $1$ and ending at the right endpoint
$1$ either agree or intersect for an elementary geometric reason. In the latter case,
splice/augment their
first crossing from the left, $d_1$. In either case, remove the path $s_1$
connecting $1$ to $1$ from the braid. (If splicing was necessary to create $s_1$, then leave
a marker $1$ on the lower strand as shown in Figure \ref{fig:marker}.) Proceed by induction
to find the
paths $s_2,\ldots,s_q$ and for those $s_i$ that were the result of splicing, leave a marker
$i$ and place the spliced crossing $d_i$ in $X$.


\begin{figure}
   \centering
   \begin{minipage}[c]{.3\textwidth}
   \centering
   \includegraphics[width=\textwidth]{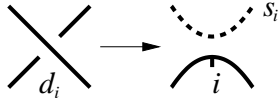}
   \end{minipage}
   \hfill
   \begin{minipage}[c]{.6\textwidth}
   \centering
   \caption{Splicing a crossing to create the path $s_i$; after removing $s_i$,
   a marker $i$ is left on the remaining diagram.}
   \label{fig:marker}
   \end{minipage}
\end{figure}

The components of $L_\beta$ are enumerated by the cycles in the permutation $\pi$ that
underlies $\beta$ and the construction treats these components independently of one another.
The number of elements in $X$ is $q$ minus the number $\mu$ of these cycles/components: it
is exactly the
largest element $i$ of each cycle of $\pi$ whose corresponding path $s_i$ `exists
automatically,' without splicing. A way to see this is the following. Suppose $\pi$ contains
a single cycle. Unless $q=1$, $d_1$ exists. When $s_1$ is removed from the braid, $1$ is
`cut out' of $\pi$: in the next, smaller braid, the underlying permutation takes
$\pi^{-1}(1)$ to $\pi(1)$. In particular, we still have a single cycle. Unless $2$ is its
largest (and only) element, $d_2$ will exist and the removal of $s_2$ cuts $2$ out of the
permutation. This goes on until we reach $q$, at which stage the braid is a single strand
and more splicing is neither possible nor necessary.

We define an oriented graph $G_\beta$ on the vertex set $\{\,1,2,\ldots,q\,\}$ by the rule
that an oriented edge connects $i$ to $j$ if $s_j$ contains the marker $i$. Note that $i<j$
is necessary for this and that each $i$ can be the starting vertex of at most one edge. For
that reason, $G_\beta$ doesn't even have unoriented cycles (consider the smallest number in
a supposed cycle). Thus, $G_\beta$ is a $\mu$--component forest. The largest element of each
tree is its only sink.

$X$ is a ruling with the $i$th eye partially bounded by the path $s_i$. These are easily
seen to satisfy Definition \ref{def:rul}: if the $i$th and $j$th eyes meet at the switch
$d_i$, then an edge connects $i$ to $j$ in $G_\beta$, thus $i<j$ and we see that in the
vertical slice through $d_i$, we have the second of the admissible configurations of Figure
\ref{fig:switch}. The value of $\theta$ for this ruling is $\mu$.


To prove that $X$ is also an augmentation, we'll check it directly using
the analysis of admissible discs in $\gamma_\beta$ from p.\
2056 of \cite{en}.
Note that each of $a_1,\ldots,a_q$ (Figure
\ref{fig:lagpic}) has a trivial admissible disc contributing $1$ to its
differential, so it suffices to show that for each $j$, there is exactly
one more admissible disc with positive corner at $a_j$ and all negative
corners at crossings in $X$. In fact we will use induction to prove the following:
\begin{itemize}
\item For each $j$, this second disc $\Pi_j$
will have either no negative corner or, if $d_j$ exists, then exactly one
negative corner at $d_j$.
\item In the admissible sequence corresponding to $\Pi_j$, $i$ appears if and only
if $G_\beta$  contains an oriented path from $i$ to $j$, and each such $i$ shows up
exactly once.
\end{itemize}

The path $s_1$ completes the boundary of an admissible disc with positive
corner at $a_1$. Because $s_1$ is removed in the first stage, no crossing
along $s_1$ other than $d_1$ will be in $X$.

Now, assume that for each $j<n$, a unique disc $\Pi_j$ exists with the
said properties. Building a non-trivial admissible disc with positive corner at $a_n$,
we start along the path $s_n$. (We will concentrate on the boundary of the admissible
disc. Proposition 6.4 of \cite{en} classifies, in terms of admissible sequences, which
of the possible paths correspond to admissible discs.) When we reach a marker $j$, we
are forced to enter $\partial\Pi_j$. Then by the inductive hypothesis, we have no other
choice but to follow $\partial\Pi_j$ until we reach $a_j$. There, we travel around the
$j$th trivial disc and continue along $\partial\Pi_j$, back to $d_j$ and $s_n$. By the
hypothesis, each $a_i$ is visited at most once, so their sequence is admissible. At the
next marker along $s_n$, a similar thing happens but using another, disjoint branch of
$G_\beta$, so the sequence stays admissible.

If $d_n$ exists, then upon reaching it, we seemingly get a choice of turning or not. If
we do turn, i.e.\ continue along $s_n$, then after a few more markers, we successfully
complete the construction of $\Pi_n$. Because all markers along $s_n$ were visited, it
has both of the required properties.

We still have to rule out the option of not turning at $d_n$. Suppose that's what we do.
Then we end up on a path $s_m$, where $m$ is the endpoint of the edge of $G_\beta$ starting
at $n$; in particular, $m>n$. We may encounter markers along $s_m$, but the previous
analysis applies to them and eventually we always return to $s_m$ and exit the braid at
the right endpoint $m$ (or at an even higher number, in case we left $s_m$ at $d_m$). But
this is impossible by Lemma 6.2 of \cite{en}.
\end{proof}

\begin{figure}
\centering
\includegraphics[width=\linewidth]{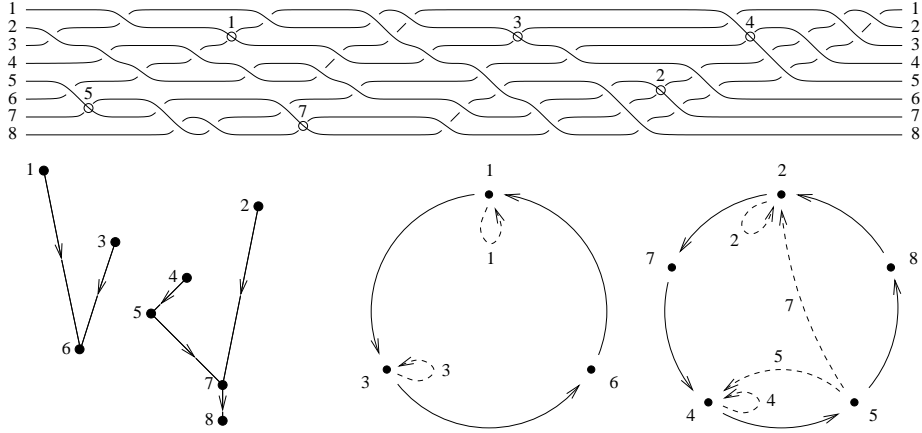}
\caption{A braid $\beta$ with an augmentation $X$ (marked crossings) which is also a ruling.
The forest graph is $G_\beta$ and the other two graph components constitute the
`graph realized by $X$,' as in \cite{en}.}
\label{fig:fonat}
\end{figure}

\begin{megj}
In \cite{en}, we used a two-component link of the braid-positive knots $8_{21}$ and
$16n_{184868}$ to illustrate a different construction of an augmentation.
Comparing Figure \ref{fig:fonat} to Figure 15 of \cite{en}, we see that the set $X$
constructed in the above proof is indeed different from that of Proposition 7.11 of
that paper. Also, the graph realized by this `new' $X$ (in the sense of Definition
7.9 in \cite{en}) is different from what we called the augmented graph of the underlying
permutation of $\beta$ there. In the example, these are both due to the fact that the
position of the augmented crossing `$3$' has changed.
\end{megj}

\section{Examples}\label{sec:ex}

The following proposition is easy to prove either using skein relations of Homfly and/or
Kauffman polynomials, or by a straightforward induction proof:


\begin{all}\label{pro:ketszal}
The ruling polynomial of the $(p,2)$ torus link is $R(z)=$
\[z^{p-1}+(p-1)z^{p-3}+{p-2\choose 2}z^{p-5}+{p-3\choose 3}z^{p-7}
+\ldots+{p-\lfloor p/2\rfloor\choose\lfloor p/2\rfloor}z^{p-2\lfloor p/2\rfloor -1}.\]
The total number of rulings is $R(1)=f_p$, the $p$'th Fibonacci number. The total number
of augmentations is $R(2^{-1/2})2^{(\chi^*+1)/2}=(2^{p+1}-(-1)^{p+1})/3$.
\end{all}

\begin{figure}
   \centering
   \begin{minipage}[c]{.5\textwidth}
   \centering
   \includegraphics[width=\textwidth]{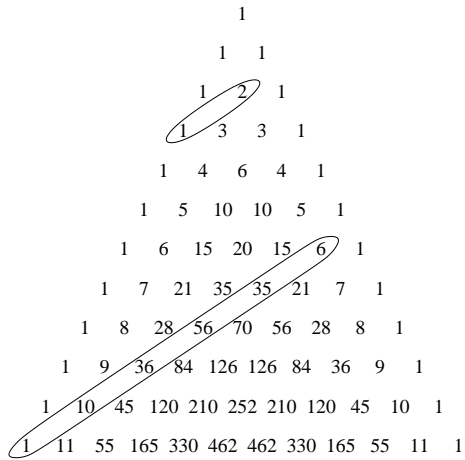}
   \end{minipage}
   \hfill
   \begin{minipage}[c]{.4\textwidth}
   \centering
   \caption{Ruling invariants of the $(3,2)$ and $(11,2)$ torus knots in Pascal's triangle.}
   \label{fig:pascal}
   \end{minipage}
\end{figure}

In particular, these ruling polynomials can be easily read off of Pascal's triangle, as
shown in Figure \ref{fig:pascal}. For example for $p=11$, we get the ruling polynomial
$R(z)=z^{10}+10z^8+36z^6+56z^4+35z^2+6$.
It seems likely that among Legendrian closures of positive braids with a given value
of $tb$, the $(p,2)$ torus link with $p=tb+2$ has the least number of rulings for all
values of $\theta$. For $tb=9$, the braid-positive knots with the largest number of
rulings (for each $\theta$) are the mutants
$13n_{981}$ and
$13n_{1104}$. These have $R(z)=z^{10}+10z^8+36z^6+60z^4+47z^2+14$.



Mutant knots share the same Kauffman and Homfly polynomials, thus mutant braid-positive
knots cannot be distinguished by their ruling polynomials.
The braid-positive knots
$12n_{679}$ and
$13n_{1176}$ are not mutants yet they share the same ruling polynomial
$R(z)=z^{10}+10z^8+36z^6+58z^4+42z^2+11$ (their Kauffman and Homfly polynomials are
actually different, but they agree in the coefficients that mean numbers of rulings).

Proposition \ref{pro:ketszal} shows that for the $(p,2)$ torus link, roughly two thirds of
the $2^p$ subsets of its crossings are augmentations. This ratio depends above all on the
number of strands in the braid and goes down approximately by a factor of two every time
the latter increases by one. When the number of strands is low, the ratio is quite
significant\footnote{Thus the relatively complicated nature of the proof of Theorem \ref{thm:szimultan} and of the construction in section 7 of \cite{en} is somewhat misleading.}. This phenomenon seems to be unique to braid-positive links. (It may be
worthwhile to compare to Chekanov's $5_2$ diagrams, where out of the $64$ subsets, only $3$,
respectively $2$, are augmentations.)

\begin{pelda}
The following were computed using a computer program written by Supap Kirtsaeng, based on
Theorem \ref{thm:gauss}. (Note that mere numbers of augmentations can also be determined
from the Homfly or Kauffman polynomials using formula \eqref{eq:aug}.) The braid word
$(\sigma_1\sigma_2)^6$, corresponding to the $(3,6)$ torus link, yields 1597 augmentations
(about $39$\% of all subsets of its crossings). The knot 
$12n_{679}$ (braid word $\sigma_1^3\sigma_2^2\sigma_1^2\sigma_2^5$) has $1653$ augmentations
(appr.\ $40$\%).
$13n_{1176}$ also has $1653$ augmentations, but its braid index is $4$; for the braid word
$\sigma_1\sigma_2^2\sigma_3\sigma_1^2\sigma_2^2\sigma_3^2\sigma_2\sigma_1\sigma_2$, the
augmentations account for only $20$\% of all subsets of crossings. The knots
$13n_{981}$ (closure of $\sigma_1\sigma_2^3\sigma_3\sigma_1\sigma_3\sigma_2^3\sigma_3^3$) and
$13n_{1104}$ ($\sigma_1\sigma_2^2\sigma_3\sigma_1\sigma_3\sigma_1^2\sigma_2^3\sigma_3\sigma_1$)
both have $1845$ augmentations (i.e., $23$\% of all possibilities work). About the following
two knots, Stoimenow \cite{stoi} found that their braid index is $4$, but in order to obtain
them as closures of positive braids, we need $5$ strands.
$16n_{92582}$ (braid word
$\sigma_1\sigma_2^2\sigma_3\sigma_4\sigma_3\sigma_1^2\sigma_2^2\sigma_3^2
\sigma_2\sigma_4\sigma_3^2$) has 7269 augmentations, which is only about $11$\% of all
possibilities.
$16n_{29507}$
($\sigma_1\sigma_2^2\sigma_3\sigma_1\sigma_3\sigma_4\sigma_1\sigma_2\sigma_4\sigma_2
\sigma_3^3\sigma_4\sigma_2$) has $8109$ ($12$\%).
\end{pelda}

\end{document}